\documentclass[11pt]{article}
\usepackage[centertags]{amsmath}
\usepackage{amsfonts}
\usepackage{amssymb}
\usepackage{amsthm}
\usepackage{soul}
\usepackage{mathrsfs}
\usepackage{extarrows}
\usepackage[margin=1in]{geometry}
\usepackage{multirow}
\usepackage{array}
\usepackage{tocloft}
\usepackage{tabularx}
\usepackage{booktabs}
\usepackage{enumitem}
\usepackage{tikz-cd}

\theoremstyle{plain}
\newtheorem{thm}{Theorem}[section]
\newtheorem{defn}[thm]{Definition}
\newtheorem{cor}[thm]{Corollary}
\newtheorem{lem}[thm]{Lemma}
\newtheorem{prop}[thm]{Proposition}
\newtheorem{exam}[thm]{Example}
\newtheorem{rem}[thm]{Remark}

\newcommand{\BB}{\mathbb{B}}
\newcommand{\CC}{\mathbb{C}}
\newcommand{\DD}{\mathbb{D}}

\newcommand{\NN}{\mathbb{N}}

\newcommand{\RR}{\mathbb{R}}

\newcommand{\BBB}{\mathcal{B}}

\newcommand{\EEE}{\mathcal{E}}

\newcommand{\HHH}{\mathcal{H}}
\newcommand{\III}{\mathcal{I}}

\newcommand{\KKK}{\mathcal{K}}

\newcommand{\PPP}{\mathcal{P}}
\newcommand{\QQQ}{\mathcal{Q}}

\newcommand{\la}{\langle}
\newcommand{\ra}{\rangle}
\newcommand{\intd}{\mathrm{d}}

\newcommand{\bz}{\bar{z}}
\newcommand{\bw}{\bar{w}}

\newcommand{\bpartial}{\bar{\partial}}
\newcommand{\cinfty}{\mathscr{C}^\infty}
\newcommand{\polytwod}{\mathbb{C}[\mathbf{z},\mathbf{\bar{z}}]}

\newcommand{\polytwo}{\mathbb{C}[z,\bar{z}]}

\newcommand{\polyzd}{\mathbb{C}[z_1,\cdots,z_d]}
\newcommand{\Tbf}{\mathbf{T}}
\newcommand{\Jbf}{\mathbf{J}}
\newcommand{\Ibf}{\mathbf{I}}
\newcommand{\Nbf}{\mathbf{N}}
\newcommand{\Sbf}{\mathbf{S}}

\newcommand{\mbf}{\mathbf{m}}
\newcommand{\Mzbf}{\mathbf{M_z}}

\newcommand{\Hol}{\mathrm{Hol}}

\allowdisplaybreaks[4]

\makeatletter\@addtoreset{equation}{section} \makeatother
\title {Operator Model via Compact Embedding}
\author{
Dexie Lin
    \thanks{College of Mathematics and Statistics, Chongqing University, 401331, Chongqing,  China, dexielin@cqu.edu.cn},
Yi Wang
	\thanks{College of Mathematics and Statistics, Center of Mathematics, Chongqing University, 401331, Chongqing,  China, wang\_yi@cqu.edu.cn}
    }
\date{ }

\setlist[enumerate,1]{label=(\arabic*)}
 \setlist[enumerate,2]{label=(\alph*)}
 \setlist[enumerate,3]{label=(\roman*)}

\begin{document}
\maketitle

\begin{abstract}
In this paper, we extend several function-theoretic and geometric constructions to the realm of multi-variable operator theory. A commuting tuple $\mathbf{T}=(T_1,\cdots,T_d)$ of bounded linear operators on a Hilbert space, equipped with a cyclic vector $h$, is abbreviated as a cyclic commuting tuple $(\mathbf{T}, h)$. We encode the complete information of such a tuple into a single positive compact operator on the Fock space.
Drawing parallels with spectral geometry, we investigate how the spectral data of this positive compact operator—its eigenvalues and eigenfunctions—reflect fundamental properties of the operator tuple.
Our main contributions are threefold. First, we establish a Weyl-type approximation formula for certain operator tuples, demonstrating that the asymptotic behavior of eigenvalues carries geometric information about the vanishing variety. Second, we construct two kernel functions derived from the eigenvalues and eigenfunctions: the first extends the Bergman kernel, while the second extends the Fourier-Laplace transformation. We prove that the Fourier-Laplace kernel defines a reproducing kernel Hilbert space (RKHS) on which the coordinate differential operators are unitarily equivalent to the adjoint tuple $\mathbf{T}^*$. Consequently, we show that the convergence points of the Bergman-type kernel characterize the joint eigenvalues of the adjoint tuple. Finally, we obtain a Paley-Wiener-Schwartz type theorem for cyclic commuting tuples, characterizing cyclic commuting tuples whose associated Agler's linear functional are distributions. For tuples consisting of matrices, we obtains a more explicit characterization.

~

	\end{abstract}
	

\section{Introduction}\label{sec: introduction}

Compact embeddings represent a fundamental concept that naturally appear throughout various areas of mathematics. A classic example arises in Riemannian geometry(cf. \cite{CrPuRa01SpectralGeometry}\cite{Kac66}): for a compact Riemannian manifold $M$, the Sobolev embedding $\iota: W^{1,2}(M) \hookrightarrow L^2(M)$ is compact. The composition $\iota\iota^\ast$ yields the Green operator, which coincides with the inverse of the Laplace--Beltrami operator (modulo the kernel of constants). The spectral data of this operator---its eigenvalues and eigenvectors---encodes deep geometric information about the underlying manifold $M$. For example, the Weyl's approximation theorem shows that the asymptotic behavior of the eigenvalues carries geometric information of the underlying manifold. Together, the eigenvalues and eigenvectors form several important kernel functions, such as the heat kernel. In this paper, we extend this fundamental principle to the realm of multi-variable operator theory. 

Assume $\Tbf=(T_1,\cdots,T_d)$ is a commuting $d$-tuple of bounded linear operators on a Hilbert space $\HHH$, i.e., $T_iT_j=T_jT_i, \forall i, j=1,\cdots,d$. Such commuting tuples are classic subject of study in multivariable operator theory. Let us assume further that $\Tbf$ has a cyclic vector $h\in\HHH$. That is, the vectors
\[
\{p(\Tbf)h~:~p\in\polyzd\}
\]
is dense in $\HHH$. For convenience, let us call $(\Tbf,h)$ a cyclic commuting ($d$-)tuple. 

Denote $F^2$ the Fock space, or the Segal-Bargmann space on $\CC^d$.
Then it is elementary to check that the the operator
\[
\iota: F^2\to\HHH,\quad f\mapsto f(\Tbf)h
\]
is compact. The positive compact operator\footnote{Here we consider $\iota^*\iota$ instead of $\iota\iota^*$, because it is on the more tractable space $F^2$. Also technically, $\iota$ may not be injective. But we keep the terminology ``compact embedding'' in the title just for convenience.} $L_{\Tbf,h}=\iota^\ast\iota$ encodes the cyclic tuple $(\Tbf,h)$ into a single operator on $F^2$. Let $L_{\Tbf,h}=\sum_k \lambda_kf_k\otimes f_k$ be its spectral decomposition. Here, $\{\lambda_k\}$ is the sequence of nonzero eigenvalues, listed in non-increasing order, also known as the singular values, and $\{f_k\}$ is a corresponding sequence of orthonormal eigenvectors. Analogous to the situation in spectral geometry, we investigate how the asymptotic behavior of 
$\{\lambda_k\}$ and $\{f_k\}$ reflect the properties of $(\Tbf,h)$.
Consider the following simple example.

\begin{exam}\label{exam: intro homogeneous quotient module}
Let $I$ be a homogeneous ideal of the polynomial ring $\polyzd$. Let $\PPP_I$ be its closure in the Drury-Arveson space $H_d^2$. Then $\PPP_I$ is a Hilbert submodule of $H_d^2$, which means $\PPP_I$ is invariant under coordinate multiplications
\[
M_{z_i}: H_d^2\to H_d^2,\quad f\mapsto z_if,\quad i=1, \cdots,d.
\]
Let $\QQQ_I=\PPP_I^\perp\cong H_d^2/\PPP_I$ be the associated Hilbert quotient module. Then $\QQQ_I$ is invariant under $M_{z_i}^*$. Denote $Q$ the orthogonal projection onto $\QQQ_I$. Consider 
\[
S_i=QM_{z_i}|_{\QQQ_I},\quad i=1,\cdots,d,\quad\text{and}\quad e=Q1.
\]
Then $(\Sbf,e)$ form a cyclic commuting $d$-tuple. Such tuples are of special importance in multivariable dilation theory. See, for example, literature on the Arveson-Douglas Conjecture \cite{Arveson98}\cite{Arveson2000}\cite{Arveson02}\cite{Arveson05}\cite{Douglas06}\cite{DoTaYu16}\cite{EE15}\cite{GuWaADsurvey}. 
\end{exam}

Under the above setting, we show that the sequence $\{\lambda_k\}$ indeed carries some geometric information.
\begin{prop}\label{prop: intro homogeneous quotient module}
Under the setting of Example \ref{exam: intro homogeneous quotient module}, let $\{\mu_k\}$ be the sequence of singular values of $L_{\Sbf,e}$. Then we have
\[
\lim_{k\to\infty}\frac{\ln k}{\ln\ln(1/\mu_k)}=\dim_\CC Z(I),
\]
where 
\[
Z(I)=\{z\in\CC^d~:~p(z)=0, \forall p\in I\}.
\]
\end{prop}
For general cyclic commuting tuple $(\Tbf,h)$, define
\[
\III_\Tbf=\{f\in F^2~:~f(\Tbf)=0\}, \quad Z(\III_\Tbf)=\{z\in\CC^d~:~f(z)=0, \forall f\in\III_\Tbf\}.
\]
If $\III_\Tbf$ is graded, i.e., $\III_\Tbf$ is the closure of a homogeneous ideal, then we show that one side of the inequality holds.
\begin{thm}\label{thm: intro Weyl type inequality}
Assume $(\Tbf,h)$ is a cyclic commuting tuple such that $\III_\Tbf$ is graded. Let $\{\lambda_k\}$ be the sequence of singular values of $L_{\Tbf,h}$. Then
\[
\limsup_{k\to\infty}\frac{\ln k}{\ln\ln(1/\lambda_k)}\leq \dim_\CC Z(\III_\Tbf).
\]
\end{thm}
We are also interested in the sequence of eigenfunctions $\{f_k\}$. To illustrate its importance, let us consider another simple example.

\begin{exam}\label{exam: intro Bergman space}
Suppose $\Omega$ is a bounded open set in $\CC^d$. Let $L_a^2(\Omega)$ be the Bergman space on $\Omega$. $L_a^2(\Omega)$ consists of all square integrable holomorphic functions on $\Omega$. Consider 
\[
T_i=M_{z_i}: L_a^2(\Omega)\to L_a^2(\Omega),~f\mapsto z_i f.
\]
Assuming that $\Omega$ is relatively nice, in the sense that $\polyzd$ is dense in $L_a^2(\Omega)$, then $(\Tbf, 1)$ is a cyclic commuting tuple. The operator $L_{\Tbf,1}$ is simply the Toeplitz operator defined by
\[
\la L_{\Tbf,1}f,g\ra_{F^2}=\la f, g\ra_{L_a^2(\Omega)}=\int_\Omega f\bar{g}\intd m,\quad\forall f, g\in F^2.
\]
In this case, we also show that
\[
\lim_{k\to\infty}\frac{\ln k}{\ln\ln(1/\lambda_k)}=n.
\]
Moreover, $L_{\Tbf,1}$ has dense range. Thus the eigenvectors $\{f_k\}$ form an orthonormal basis of $F^2$. Also by the equation above,
\[
\la f_j, f_k\ra_{L_a^2(\Omega)}=\la L_{\Tbf,1}f_j,f_k\ra_{F^2}=\delta_{j,k}\lambda_j.
\]
Thus $\{\frac{f_k}{\sqrt{\lambda_k}}\}$ form an orthonormal basis of $L_a^2(\Omega)$. As a consequence,
\begin{enumerate}
    \item $\sum_k\lambda_k^{-1}f_k(z)\overline{f_k(w)}=K^\Omega(z,w)$, the Bergman kernel on $\Omega$;
    \item $\sum_k f_k(z)\overline{f_k(w)}=K^{F^2}(z,w)=e^{\la z,w\ra}$, the Fock kernel.
\end{enumerate}
This shows that the functions $\{f_k\}$ play a similar role as the monomials for $\BB_d$.
We also consider
\begin{enumerate}[start=3]
    \item $F(z,w)=\sum_k\lambda_k f_k(z)\overline{f_k(w)}=\la L_{\Tbf,1}K_w^{F^2},K_z^{F^2}\ra_{F^2}$, the kernel function of $L_{\Tbf,1}$.
\end{enumerate}
Another elementary observation is the following.
\begin{enumerate}[start=4]
      \item up to a change of variable, $F(z,w)$ is the Fourier-Laplace transform (cf. \cite{HormanderPDEI2003}) of $\chi_\Omega$:
      \[
      F(z,w)=\widehat{\chi_\Omega}\left((i(\bw+z),-\bw+z)\right),\quad\forall z, w\in\CC^d.
      \]
  \end{enumerate}
\end{exam}
Having the above observations in mind, for general tuple $(\Tbf,h)$, we are interested in the following three kernel functions:
\[
\sum_k\lambda_k^{-1}f_k(z)\overline{f_k(w)},\quad\sum_kf_k(z)\overline{f_k(w)},\quad\sum_k\lambda_kf_k(z)\overline{f_k(w)}.
\]
From the expression, it is immediate that these are (semi-)positive definite kernels. Thus the theory of reproducing kernel Hilbert space \cite{PaRaRKHSbook16} may play a role. For a cyclic commuting tuple $(\Tbf,h)$, let us analogously define
\[
F_{\Tbf,h}(z,w)=\sum_k\lambda_kf_k(z)\overline{f_k(w)}=\la L_{\Tbf,h}K^{F^2}_w, K^{F^2}_z\ra.
\]
The theorem below uses $F_{\Tbf,h}$ to give an operator model for $\Tbf^\ast=(T_1,\cdots,T_d)$.

\begin{thm}\label{thm: intro H(F) and T}
Suppose $(\Tbf,h)$ is a cyclic commuting $d$-tuple. Denote $\HHH(F_{\Tbf,h})$ the reproducing kernel Hilbert space on $\CC^d$ defined by $F_{\Tbf,h}$. Then the following hold.
\begin{enumerate}
    \item $\HHH(F_{\Tbf,h})\subset F^2$ as sets of functions on $\CC^d$;
    \item the functions $\{\sqrt{\lambda_k}f_k\}$ form an orthonormal basis for $\HHH(F_{\Tbf,h})$;
    \item the differentiation operators
    \[
    \partial_i: \HHH(F_{\Tbf,h})\to\HHH(F_{\Tbf,h}),\quad f\mapsto\frac{\partial f}{\partial z_i},\quad i=1,\cdots,d
    \]
    are bounded on $\HHH(F_{\Tbf,h})$; moreover,
    \item there is a unitary operator $U: \HHH\to\HHH(F_{\Tbf,h})$ such that
    \[
    T_i^\ast=U^\ast\partial_i U,\quad i=1,\cdots,d.
    \]
\end{enumerate}
\end{thm}
The following corollary of Theorem \ref{thm: intro H(F) and T} gives a partial analogue of observation (1) in Example \ref{exam: intro Bergman space}.
\begin{cor}\label{cor: intro joint eigenvalue for T*}
Suppose $(\Tbf,h)$ is a cyclic commuting $d$-tuple. Then we have
\[
\left\{z\in Z(\III_\Tbf):\sum_k\lambda_k^{-1}|f_k(z)|^2<\infty\right\}=\overline{\sigma_p(\Tbf^*)}.
\]
\end{cor}
In light of Observation (4) in Example \ref{exam: intro Bergman space}, $F_{\Tbf,h}$ can also be viewed as an extension of the Fourier-Laplace transform to the multivariable operator setting. Recall that the Paley-Wiener-Schwartz Theorem characterizes compactly supported distributions from their Fourier-Laplace transforms. We provide an analogue of it to this setting. Denote $\polytwod=\CC[z_1,\cdots,z_d,\bar{z_1},\cdots,\bar{z_d}]$. For a cyclic commuting $d$-tuple $(\Tbf,h)$, define its associated Agler's linear functional to be the linear functional $\Lambda_{\Tbf,h}$ on $\polytwod$, determined by the moment sequence
\[
\Lambda_{\Tbf,h}(z^\alpha\bz^\beta)=\la \Tbf^\alpha h, \Tbf^\beta h\ra,\quad\alpha, \beta\in\NN_0^d.
\]
Equivalently,
\[
\Lambda_{\Tbf,h}(p\bar{q})=\la q(\Tbf)^\ast p(\Tbf)h,h\ra,\quad\forall p, q\in\polyzd.
\]
In the case $d=1$, spectral theory of normal operators and the Bram-Halmos criterion for subnormality (cf. \cite[Theorem 1.9]{ConwaySubnormal}) shows that the following are equivalent.
\begin{enumerate}
    \item[(i)] $T$ is subnormal;
    \item[(ii)] $\Lambda_{T,h}$ is a finite compactly supported positive measure $\mu$ on $\CC$, i.e., $\Lambda_{\Tbf,h}(\phi)=\int\phi\intd\mu, \forall\phi\in\polytwo$ ;
    \item[(iii)] $\Lambda_{T,h}(|\phi|^2)\geq0, \forall\phi\in\polytwo$.
    \end{enumerate} 
Note that for a (sufficiently large) compact $K$ in $\CC$, $\polytwo$ is a dense subspace of $C(K)$. Thus the measure $\mu$, if exists, is unique. 
The linear functional $\Lambda_{\Tbf,h}$ was first introduced by Jim Agler in \cite{Agler82}. Since then it has become a useful tool in the study of ``near subnormal operators''.  Here the ``near subnormality notions'' include hyponormality, $k$-hyponormality, polynomial hyponormality, etc. We refer to \cite{Agler85}\cite{Agler88}\cite{CuPuPHO}\cite{CuPu93}\cite{McPa89} for some results in this direction. It also allows connections between operator theory and other problems such as Positivstellensatz, sum of squares, moment problems and control theory (cf. \cite{HePu07}\cite{Putinar93}). 

Meanwhile, note that $\polytwod$ is also a dense subspace of $\cinfty(\CC^d)$, whose dual space consists exactly of all compactly supported (Schwartz) distributions on $\CC^d$. Thus it makes sense to ask: for what $(\Tbf,h)$ is $\Lambda_{\Tbf,h}$ a distribution? Examples include Jordan blocks, coordinate muliplication operators on the Drury-Arveson space, and the class of subjordan operators (cf. Remark \ref{rem: subjordan}). 
The following theorem can be viewed as an extension of the Paley-Wiener-Schwartz theorem to the realm of multivariable operator theory.
\begin{thm}\label{thm: intro PWS}
Suppose $(\Tbf, h)$ is a cyclic commuting $d$ tuple, $N\in\NN_0$, and $K\subset\CC^d$ is a compact convex set. Then the following are equivalent.
\begin{enumerate}
    \item $\Lambda_{\Tbf,h}\in\left(\cinfty(\CC^d)\right)'$, with order $\leq N$, supported in $K$;
    \item there exists a constant $C>0$ such that
\begin{equation}\label{eqn: intro F PWS}
|F_{\Tbf,h}(z,w)|\leq C(1+|z|+|w|)^Ne^{H_K(z+w)},\quad\forall z, w\in\CC^n.
\end{equation}
Here $H_K$ denotes the complex support function of $K$ (cf. Definition \ref{defn: complex supporting function}).
\end{enumerate}
\end{thm}
In the diagonal direction, Inequality \eqref{eqn: intro F PWS} indicates that $F_{\Tbf,h}(z,z)$ is of exponential growth, which is not surprising by the definition of $F_{\Tbf,h}$. Thus Inequality \eqref{eqn: intro F PWS} essentially says that $F_{\Tbf,h}$ is of polynomial growth in the off-diagonal direction, i.e., $|F_{\Tbf,h}(z,-z)|\leq C(1+|z|)^N$. Together with Theorem \ref{thm: intro H(F) and T}, this gives a characterization and a model of such operator tuples. In the case when $\Tbf$ consists of matrices, the characterization becomes simpler. We address this in Theorem \ref{thm: matrix case}.

This paper is organized as follows. In the Preliminary section, we fix a few terminologies and recall some basic facts. In Section \ref{sec: Weyl}, we give the proofs of Proposition \ref{prop: intro homogeneous quotient module} and Theorem \ref{thm: intro Weyl type inequality}. In Section \ref{sec: operator model and spectrum}, we give the proofs of Theorem \ref{thm: intro H(F) and T} and Corollary \ref{cor: intro joint eigenvalue for T*}. In Section \ref{sec: PWS type theorem}, we give a few examples when the Agler's linear functional is given by a distribution, and then we prove Theorem \ref{thm: intro PWS}. In Section \ref{sec: the matrix case}, we give a characterization of cyclic commuting matrix tuples given by a distribution.

The observations in this paper lead to a few further interesting questions. We list a few as example. First, like in spectral geometry, one may ask whether there exist isospectral tuples (cyclic commuting tuples with the same $\{\lambda_k\}$). Second, in light of Corollary \ref{cor: intro joint eigenvalue for T*}, can one characterize other points in the Taylor spectrum $\sigma(\Tbf)$ using asymptotic behavior of $\{\lambda_k\}$ and $\{f_k\}$? Finally, can the generalized Fourier transform $F_{\Tbf,h}$ be used to characterize the compactness, Schatten-class membership of hereditary functional calculus of $\Tbf$ (cf. Remark \ref{rem: model for hereditary})? Of course, none of these are easy questions. But we believe that addressing these questions may offer some complementary insights to the current operator model theory framework.

\noindent{\bf Acknowledgment:} The authors would like to thank J.William Helton, Catalin Badea, Ra\'{u}l Curto, Orr Shalit, Joseph Ball, Zipeng Wang, Penghui Wang, Rongwei Yang, Kui Ji, Kunyu Guo for helpful discussions. Special thanks goes to Michael Hartz and Mihai Putinar for reading a previous version of this paper carefully and providing valuable suggestions.

\section{Preliminary}\label{sec: preliminary}
In this section, we fix a few terminologies and recall some basic facts. 

\begin{defn}\label{defn: cyclic commuting tuple}
Let $\HHH$ be a Hilbert space. By a \emph{cyclic commuting ($d$-)tuple on $\HHH$}, we mean a $d$-tuple $\Tbf=(T_1,\cdots,T_d)$ of pairwise commuting bounded linear operators on $\HHH$, together with a cyclic vector $h\in\HHH$ for $\Tbf$. In other words, $(\Tbf,h)$ is a cyclic commuting $d$-tuple on $\HHH$ if 
\[
T_i\in\BBB(\HHH),~T_iT_j=T_jT_i,~i, j=1,\cdots,d,\quad\text{and}\quad\overline{\{p(\Tbf)h~:~p\in\polyzd\}}=\HHH.
\]
It will be the main subject of study in this paper. There is then a natural concept of unitary equivalence between cyclic commuting tuples. Suppose $(\Sbf,e)$ is another cyclic commuting $d$-tuple on a Hilbert space $\EEE$. We say that $(\Tbf,h)$ is unitarily equivalent to $(\Sbf,e)$, denoted $(\Tbf,h)\cong(\Sbf,e)$, if there is a unitary operator $U:\HHH\to\EEE$, such that
\[
U^*S_iU=T_i,~i=1,\cdots,d\quad\text{and}\quad Uh=e.
\]
\end{defn}

Let $F^2$ be the Fock space on $\CC^d$, which we introduce in Subsection \ref{subsec: analytic function spaces}.

\begin{defn}\label{defn: iota L the sequence and F}
Suppose $(\Tbf,h)$ is a cyclic commuting $d$-tuple on $\HHH$. Define the mapping
\[
\iota=\iota_{\Tbf,h}: F^2\to\HHH,\quad f\mapsto f(\Tbf)h,
\]
and let $L_{\Tbf,h}=\iota^*\iota$. In Lemma \ref{lem: iota is compact} below, we show that $\iota$ is compact. As a consequence, $L_{\Tbf,h}$ is a positive compact operator on $F^2$. Let
\[
L_{\Tbf,h}=\sum_k\lambda_k f_k\otimes f_k,\quad \lambda_1\geq\lambda_2\geq\cdots>0
\]
be its spectral decomposition, where $\{\lambda_k\}$ is the sequence of non-zero eigenvalues, listed in non-increasing order, and $\{f_k\}$ is an orthonormal sequence of eigenvectors. In the case of repeated eigenvalues, the choice of $\{f_k\}$ is not unique. We fix one such choice. Define the kernel function
\[
F_{\Tbf,h}(z,w)=\la L_{\Tbf,h}K^{F^2}_w, K_z^{F^2}\ra=\sum_k\lambda_kf_k(z)\overline{f_k(w)},\quad z, w\in\CC^d.
\]
From the last expression, it is clear that $F_{\Tbf,h}$ is a positive definite kernel function on $\CC^d$. (See Subsection \ref{subsec: RKHS} below.) Denote $\HHH(F_{\Tbf,h})$ the reproducing kernel Hilbert space on $\CC^d$ determined by $F_{\Tbf,h}$. Also, let us define the semi-inner product
\[
\la p, q\ra_{\Tbf,h}=\la p(\Tbf)h, q(\Tbf)h\ra_\HHH,\quad\forall p, q\in\polyzd,
\]
and the linear functional
\[
\Lambda_{\Tbf,h}:\CC[z_1,\cdots,z_d,\bar{z_1},\cdots,\bar{z_d}]\to\CC,\quad\Lambda_{\Tbf}(z^\alpha\bar{z}^\beta)=\la z^\alpha, z^\beta\ra_{\Tbf,h},\quad\forall \alpha, \beta\in\NN_0^d.
\]
It is clear that each of $L_{\Tbf,h}$, $F_{\Tbf,h}$, $\HHH(F_{\Tbf,h})$, $\Lambda_{\Tbf,h}$, and $\la\cdot,\cdot\ra_{\Tbf,h}$ determines the others. Below, we show that by a GNS-type construction, the semi-inner product $\la\cdot,\cdot\ra_{\Tbf,h}$ determines $(\Tbf,h)$ up to unitary equivalence. Thus each of the five objects contains complete information of $(\Tbf,h)$.
In the sequel, when no confusion is caused, we omit the subscripts and simply write $\iota, L, F, \HHH(F), \Lambda$.
\end{defn}

\noindent{\bf the GNS-type Construction:} Let $\la\cdot,\cdot\ra_{\Tbf,h}$ be as in Definition \ref{defn: iota L the sequence and F}. Denote $\|\cdot\|_{\Tbf,h}$ the induced semi-norm. Let
\[
I=\{p\in\polyzd:\|p\|_{\Tbf,h}=0\}=\{p\in\polyzd:\|p(\Tbf)h\|=0\}=\{p\in\polyzd:p(\Tbf)=0\}.
\]
The last equality follows from the cyclicity of $h$.
Then $I$ is an ideal of $\polyzd$. Still denote by $\la\cdot,\cdot\ra_{\Tbf,h}, \|\cdot\|_{\Tbf,h}$ the induced inner product and norm on $\polyzd/I$. The mappings
\[
M_{z_i}: \polyzd/I\to\polyzd/I, \quad p+I\mapsto z_ip+I,\quad i=1,\cdots,d
\]
are well-defined. Also define
\[
U: (\polyzd/I, \|\cdot\|_{\Tbf,h})\to\HHH,\quad p+I\mapsto p(\Tbf)h.
\]
Then it is easy to check that $U$ is an isometry, and the following commuting diagram holds.
\begin{center}
    \begin{tikzcd}
        \polyzd/I \ar[r,"M_{z_i}"] \ar[d,"U"] &\polyzd/I \ar[d,"U"]\\ \HHH \ar[r,"T_i"]&\HHH
    \end{tikzcd}
\end{center}
Let $\HHH_{\Tbf,h}$ be the completion of $(\polyzd/I, \|\cdot\|_{\Tbf,h})$. Then $M_{z_i}$ extends to a bounded linear operator on $\HHH_{\Tbf,h}$, and $(\Mzbf, 1+I)\cong(\Tbf,h)$.

\subsection{Analytic Function Spaces}\label{subsec: analytic function spaces}
In this subsection, we recall some basics about the Bergman space, the Fock space and the Drury-Arveson space. For more about these spaces, we refer the readers to \cite{ZhuFockBook12}\cite{Ha23DAspace}.
\begin{defn}
Denote $\intd m$ the Lebesgue measure on $\CC^d$. For an open set $\Omega$ in $\CC^d$, denote $\Hol(\Omega)$ the space of analytic functions on $\Omega$. The Bergman space $L_a^2(\Omega)=L^2(\Omega)\cap\Hol(\Omega)$. Denote $\intd\mu(z)=\pi^{-d}e^{-|z|^2}\intd m(z)$ the Gaussian measure on $\CC^d$. The Fock space is defined by
\[
F^2=\{f\in\Hol(\CC^d)~:~\|f\|_{F^2}^2=\int|f|^2\intd\mu<\infty\}.
\]
Finally, the Drury-Arveson space $H_d^2$ is the space of analytic functions on the open unit ball $\BB_d$ with Taylor expansion
\[
f(z)=\sum_{\alpha\in\NN_0^d}a_\alpha z^\alpha,\quad\text{where}\quad\|f\|_{H_d^2}^2=\sum_\alpha|a_\alpha|^2\frac{\alpha!}{|\alpha|!}.
\]
Temporarily denote $\HHH$ for either of the three spaces. It is well-known that $\HHH$ is a reproducing kernel Hilbert spaces on relative domains. This means for any $z\in\Omega$, $\CC^d$, or $\BB_d$, respectively, the evaluation map
\[
\mathrm{ev}_z: \HHH\to\CC,  f\mapsto f(z)
\]
defines a bounded linear functional on $\HHH$. By the Riesz representation theorem, there is an element in the space, denoted by $K_z^\HHH$, such that
\[
f(z)=\la f, K_z^\HHH\ra_\HHH,\quad\forall f\in\HHH.
\]
It is also conventional to write $K^\HHH(w,z)=K^\HHH_z(w)$. We have
\[
K^{F^2}(z,w)=e^{\la z, w\ra},~z, w\in\CC^d\quad\text{and}\quad K^{H_d^2}(z,w)=\frac{1}{1-\la z,w\ra},~z, w\in\BB_d.
\]
Suppose $\{e_n\}$ is an orthonormal basis for $\HHH$. Then a standard argument shows that
\[
K^\HHH(w,z)=\sum_ne_n(w)\overline{e_n(z)}.
\]
\end{defn}

The following lemma is elementary for experts. We provide a proof for completeness.
\begin{lem}\label{lem: iota is compact}
Let $(\Tbf,h)$ and $\iota$ be as in Definition \ref{defn: iota L the sequence and F}. Then $\iota$ is compact.
\end{lem}
\begin{proof}
Choose $r>0$ so that $\|T_i\|<\frac{r}{2}, i=1,\cdots,d$. Then there exists a constant $C_1, C_2>0$ such that
\[
\|\iota f\|=\|f(\Tbf)h\|\leq\|h\|\|f(\Tbf)\|\leq C_1\|f\|_{r\BB_d}\leq C_2 \|f\|_{L_a^2(2r\BB_d)},\quad \forall f\in F^2.
\]
Here $\|\cdot\|_{r\BB_d}$ denote the supremum norm on $r\BB_d$, and $L_a^2(2r\BB_d)$ denotes the Bergman space on $2r\BB_d$. Since $L_a^2(2r\BB_d)$ and $F^2$ shares a same orthogonal basis $\{z^\alpha\}_{\alpha\in\NN_0^d}$, comparing the norms shows that the restriction map $F^2\to L_a^2(2r\BB_d)$ is compact. Consequently, $\iota$ is also bounded.
\end{proof}

\subsection{Reproducing Kernel Hilbert Spaces}\label{subsec: RKHS}
Let us also recall some basic facts in the general theory of reproducing kernel Hilbert spaces (RKHS). Our main source of references is \cite{PaRaRKHSbook16}. Let $\Omega$ be any set. A function of two variables
\[
F: \Omega\times\Omega\to\CC
\]
is a positive definite kernel function on $\Omega$ provided that for any finite collection $\{z_i\}\subset\Omega$, the finite matrix $[F(z_i, z_j)]$ is a semi-positive definite matrix. Here the terminology is different from \cite{PaRaRKHSbook16}, where such $F$ is referred to simply as kernel function. We use the extended terminology to avoid confusion with the kernel functions of operators on analytic function spaces, which are also mentioned in this paper. Given a positive definite kernel function $F$ on $\Omega$, there is a unique reproducing kernel Hilbert space $\HHH(F)$ associated to it. That is, 
\begin{enumerate}
    \item $\HHH(F)$ consists of functions on $\Omega$;
    \item for any $z\in\Omega$, the evaluation maps
    \[
    \mathrm{ev}_z: \HHH(F)\to\CC,\quad f\mapsto f(z)
    \]
    is a bounded linear functional on $\HHH(F)$;
    \item write $F_w(z):=F(z,w)$. Then $F_w\in\HHH(F), \forall w\in\Omega$, and $F_w$ is the reproducing kernel of $\HHH(F)$ at $w$:
    \[
    f(w)=\la f, F_w\ra_{\HHH(F)},\quad\forall f\in\HHH(F), \forall w\in\Omega.
    \]
\end{enumerate}
\begin{lem}[\cite{PaRaRKHSbook16} Theorem 3.11]\label{lem: f in RKHS criterion}
Let $\HHH$ be an RKHS on a set $X$ with reproducing kernel $K$ and let $f: X\to\CC$ be a function. Then the following are equivalent:
\begin{enumerate}
    \item $f\in\HHH$;
    \item there exists a constant $c\geq0$, such that for every finite subset $F=\{x_1,\cdots,x_m\}\subset X$, there exists a function $h\in\HHH$ with $\|h\|\leq c$ and $f(x_i)=h(x_i), i=1,\cdots,m$;
    \item there exists a constant $c\geq0$ such that the function $c^2K(x,y)-f(x)\overline{f(y)}$ is a positive definite kernel function.
\end{enumerate}
Moreover, if $f\in\HHH$ then $\|f\|$ is the least $c$ that satisfies the inequalities in (2) and (3).
\end{lem}

\begin{defn}
Let $\HHH$ be a Hilbert space with inner product $\la\cdot,\cdot\ra$. A set of vectors $\{f_s\}_{s\in S}$ is called a \emph{Parseval frame} for $\HHH$ provided that 
\[
\|h\|^2=\sum_{s\in S}|\la h, f_s\ra|^2.
\]
\end{defn}

\begin{lem}[\cite{PaRaRKHSbook16} Proposition 2.9]\label{lem: Parseval frame criterion}
Let $\{f_s\}_{s\in S}$ be a Parseval frame for a Hilbert space $\HHH$, then there is a Hilbert space $\KKK$ containing $\HHH$ as a subspace and an orthonormal basis $\{e_s\}_{s\in S}$ for $\KKK$, such that $f_s=P_\HHH(e_s), s\in S$, where $P_\HHH$ denotes the orthogonal projection of $K$ onto $\HHH$.
\end{lem}

\begin{lem}[\cite{PaRaRKHSbook16} Theorem 2.10]\label{lem: Parseval frame in RKHS}
Let $\HHH$ be an RKHS on $X$ with reproducing kernel $K$ and let $\{f_s\}_{s\in S}\subset\HHH$. Then $\{f_s\}_{s\in S}$ is a Parseval frame for $\HHH$ if and only if $K(x,y)=\sum_{s\in S}f_s(x)\overline{f_s(y)}$, where the series converges pointwise.
\end{lem}

From Lemmas \ref{lem: Parseval frame criterion} and \ref{lem: Parseval frame in RKHS}, we have the following criterion of orthonormal basis of an RKHS.
\begin{lem}\label{lem: ONB in RKHS}
Let $\HHH$ be an RKHS on $X$ with reproducing kernel $K$ and let $\{f_s\}_{s\in S}\subset\HHH$. Then $\{f_s\}_{s\in S}$ form an orthonormal basis for $\HHH$ if and only if 
\[
\|f_s\|=1,~\forall s\in S,\quad\text{and}\quad K(x,y)=\sum_{s\in S}f_s(x)\overline{f_s(y)}.
\]
\end{lem}
\begin{proof}
By the second condition, $\{f_s\}_{s\in S}$ form an Parseval frame for $\HHH$. Let $\KKK, \{e_s\}_{s\in S}$ be as in Lemma \ref{lem: Parseval frame criterion}. By the first condition, $f_s=e_s$. Thus $\{f_s\}_{s\in S}$ form an orthonormal basis. This completes the proof.
\end{proof}

\section{Weyl Type Estimates}\label{sec: Weyl}
In this section, we give the proofs of Proposition \ref{prop: intro homogeneous quotient module} and Theorem \ref{thm: intro Weyl type inequality}. 
\begin{proof}[{\bf Proof of Proposition \ref{prop: intro homogeneous quotient module}}]
Let $I$, $\PPP_I$, $\QQQ_I$, $\Sbf, e, \{\mu_k\}, p$ be as in Example \ref{exam: intro homogeneous quotient module} and Proposition \ref{prop: intro homogeneous quotient module}. Denote $p=\dim_\CC Z(I)$. Since $I$ is homogeneous, $\QQQ_I$ has an orthogonal decomposition $\QQQ_I=\oplus_n \QQQ_n$, where $\QQQ_n$ consists of homogeneous polynomials of degree $n$. Observe that $\|z^\alpha\|_{H_d^2}^2=\frac{\alpha!}{|\alpha|!}=\frac{\|z^\alpha\|_{F^2}^2}{|\alpha|!}$. So for a homogeneous polynomial $q$ of degree $n$, one has $\|q\|_{H_d^2}^2=\frac{\|q\|_{F^2}^2}{n!}$. For $f\in F^2$ with homogeneous decomposition $f=\sum_nf_n$,
\[
\la L_{\Sbf,e}f,f\ra_{F^2}=\|\iota f\|_{\QQQ_I}^2=\|Qf\|_{H_d^2}^2=\sum_n\|Q_nf\|_{H_d^2}^2=\sum_n\frac{1}{n!}\|Q_n f\|_{F^2}^2.
\]
Here we abuse notation and denote $Q_n$ for both the orthogonal projection from $F^2$ onto $\QQQ_n$ and from $H_d^2$ onto $\QQQ_n$. Note that they coincide on $F^2$. From the above,
\[
L_{\Sbf,e}=\sum_{n=0}^\infty\frac{1}{n!}Q_n.
\]
Denote $H(n)=\dim \QQQ_n$. Then
\[
\mu_k=\frac{1}{n!},\quad\text{where }H(0)+\cdots+H(n-1)<k\leq H(0)+\cdots+H(n).
\]
On the other hand, the Hilbert-Serre theorem \cite[Chapter I, Theorem 7.5]{Hartshorne} ensures a polynomial $h$ of degree $p-1$, such that for $n$ large enough, $H(n)=h(n)$. Consequently,
\[
H(0)+\cdots+H(n)\approx n^p.
\]
Therefore
\[
\frac{\ln\left(H(0)+\cdots+H(n-1)\right)}{\ln\ln n!}<\frac{\ln k}{\ln\ln(1/\mu_k)}\leq\frac{\ln\left(H(0)+\cdots+H(n)\right)}{\ln\ln n!}.
\]
As $k\to\infty$, both sides converge to $p$. This completes the proof.
\end{proof}

\begin{proof}[{\bf Proof of Theorem \ref{thm: intro Weyl type inequality}}]
Let $(\Tbf,h)$ be any cyclic commuting $d$-tuple such that
\[
\III=\{f\in F^2~:~f(\Tbf)=0\}
\]
is graded. Let $I=\III\cap\polyzd$. Then $I$ is a homogeneous ideal and $Z(I)=Z(\III)$. Denote $(\Sbf, e)$ the tuple induced by $I$ as in Example \ref{exam: intro homogeneous quotient module}. Let $Q_n, \{\mu_k\}$ be as in the proof of Proposition \ref{prop: intro homogeneous quotient module}. Choose $r>0$ so that $r\|T_i\|<\frac{1}{2}, \forall i$. For any $f\in\polyzd$, denote $f_{1/r}(z)=f(z/r)$. Then
\[
\la L_{\Tbf,h}f,f\ra_{F^2}=\|f(\Tbf)h\|^2=\|f_{1/r}(r\Tbf)h\|^2\leq C\|f_{1/r}\|_{H^2_d}^2,
\]
where $C>0$ is a constant. Replacing $f$ with $f+g$ for any $g\in I$ and taking infimum gives
\[
\la L_{\Tbf,h}f,f\ra_{F^2}\leq C\|Q_I f_{1/r}\|^2_{H_d^2}=C\sum_{n=0}^\infty\frac{1}{r^{2n}n!}\|Q_nf\|_{F^2}.
\]
From this we see that 
\[
\lambda_k\leq \frac{C}{r^{2n}n!},\quad\text{where }H(0)+\cdots+H(n-1)<k\leq H(0)+\cdots+H(n).
\]
Continuing as in the proof of Proposition \ref{prop: intro homogeneous quotient module} gives the desired inequality.
\end{proof}

In the case of Example \ref{exam: intro Bergman space}, i.e., when $\HHH$ is the Bergman space and $\Tbf$ is the tuple of coordinate multiplication operators, then $L$ is a Toeplitz operator on the Fock space. The asymptotic behavior of Toeplitz operators on Fock and Fock type spaces was studied in \cite{ElMaMaNa16}\cite{HuWaZe26}.

\begin{prop}\label{prop: Toeplitz asymptotic}
Under the setting of Example \ref{exam: intro Bergman space}, we have
\[
\lim_{k\to\infty}\frac{\ln k}{\ln\ln(1/\lambda_k)}=n.
\]
\end{prop}
\begin{proof}
Let us denote $T_\phi$ for the Toeplitz operator
\[
T_\phi: F^2\to F^2,\quad \la T_\phi f, g\ra=\int\phi(z)f(z)\overline{g(z)}\intd m(z).
\]
Assume first that $\Omega=r\BB_d$ for some $r>0$. Then $L=T_{\chi_{r\BB_d}}$.
\[
\la Lz^\alpha,z^\beta\ra=\int_{r\BB_d}z^\alpha\bz^\beta\intd m(z)=\delta_{\alpha,\beta}r^{2|\alpha|+2d}\frac{c\alpha!d!}{(|\alpha|+d)!},
\]
where $c>0$ is a normalizing constant. Denote $P_n$ the space of homogeneous polynomials of degree $n$. We also use $P_n$ to denote the orthogonal projection onto $P_n$.
By a similar proof of Proposition \ref{prop: intro homogeneous quotient module}, we have 
\[
L=\sum_{n=0}^\infty\frac{cr^{2n+2d}d!}{(n+d)!}P_n.
\]
Let $H(n)=\dim P_n={d+n-1\choose n-1}$. Then
\[
\lambda_k=\frac{cr^{2n+2d}d!}{(n+d)!}\quad\text{for}\quad H(0)+\cdots+H(n-1)<k\leq H(0)+\cdots+H(n).
\]
Continuing as in the proof of Proposition \ref{prop: intro homogeneous quotient module} gives the equation for $\Omega=r\BB_d$.

If $0\in\Omega$, then there exists $0<r<R$ such that $T_{\chi_{r\BB_d}}\leq L\leq T_{\chi_{r\BB_d}}$. Then the equation follows. 

In the general case, take $z\in\Omega$. Recall that 
\[
U: F^2\to F^2,\quad f(w)\mapsto f(w-z)k_z(w)
\]
defines a unitary operator, and $U^*T_{\chi_\Omega}U=T_\phi$, where
\[
\phi(\lambda)=\chi_{\Omega-z}(\lambda)e^{-|\lambda|^2+|\lambda+z|^2}.
\]
Thus there exists $0<c<C, 0<r<R$, such that 
\[
cT_{r\BB_d}\leq U^*T_{\chi_\Omega}U\leq CT_{R\BB_d}.
\]
From this it is easy to see that the equation holds. This completes the proof.
\end{proof}

\begin{rem}\label{rem: model for hereditary}
There is another motivation for studying the asymptotic behavior of $\{\lambda_k\}$. Let $h=\sum a_{\alpha,\beta}z^\alpha\bw^\beta$ be a hereditary function \cite{AgMcPickInterp}. For convenience, we assume that $h$ is a finite sum. The hereditary functional calculus of $\Tbf$ is 
\[
h(\Tbf,\Tbf^\ast)=\sum a_{\alpha,\beta}\left(\Tbf^*\right)^\beta\Tbf^\alpha.
\]
For any $f, g\in F^2$ such that $f\polyzd\subset F^2, g\polyzd\subset F^2$, 
\begin{flalign*}
\la h(\Tbf,\Tbf^*)\iota f, \iota g\ra_\HHH=&\sum a_{\alpha,\beta}\la\iota(z^\alpha f), \iota(z^\beta g)\ra=\sum a_{\alpha,\beta}\la Lz^\alpha f, z^\beta g\ra\\
=&\int_{\CC^d}\int_{\CC^d}\sum a_{\alpha,\beta}z^\alpha f(z)F(w,z)\bw^\beta \overline{g(w)}\intd\mu(z)\intd\mu(w)\\
=&\int_{\CC^d}\int_{\CC^d} h(z,w)F(w,z)f(z)\overline{g(w)}\intd\mu(z)\intd\mu(w).
\end{flalign*}
In Lemma \ref{lem: ONBs} and Lemma \ref{lem: f_j closed under polynomial} below we show that $f_k\polyzd\subset F^2, \forall k$, and $\{\lambda_k^{-1/2}\iota f_k\}$ form an orthonormal basis of $\HHH$. Therefore under this basis, $h(\Tbf,\Tbf^*)$ has matrix expression $[\lambda_j^{-1/2}\lambda_k^{-1/2}a_{j,k}]_{j,k}$, where
\[
a_{j,k}=\la h(\Tbf,\Tbf^*)\iota f_j, \iota f_k\ra=\int_{\CC^d}\int_{\CC^d}h(z,w)F(w,z)f_j(z)\overline{f_k(w)}\intd\mu(z)\intd\mu(w).
\]
In other words, $h(\Tbf,\Tbf^*)$ is the composition of an integral operator and a Schur multiplier. Thus explicit estimates of $\{\lambda_k\}$ and the kernel $F(w,z)$ can directly apply to the study of hereditary functional calculus. Indeed, a similar idea was used in \cite[Section 4]{Howe80} to study Toeplitz operators on the Bergman space on the unit ball.
\end{rem}

\section{Operator Model and Spectrum}\label{sec: operator model and spectrum}
In this section, we give the proofs of Theorem \ref{thm: intro H(F) and T} and Corollary \ref{cor: intro joint eigenvalue for T*}.

\begin{lem}\label{lem: ONBs}
Let $(\Tbf,h)$ be a cyclic commuting $d$-tuple on a Hilbert space $\HHH$, and let $\iota, L, \{\lambda_k\}, \{f_k\}, F, \HHH(F)$ be as in Definition \ref{defn: iota L the sequence and F}. Let $\III=\III_\Tbf$ be as in the introduction. Then the following hold.
\begin{enumerate}
    \item The sequence $\{\lambda_k^{-\frac{1}{2}}\iota f_k\}$ form an orthonormal basis of $\HHH$;
    \item the sequence $\{f_k\}$ form an orthonormal basis for $F^2\ominus\III$;
    \item the sequence $\{\sqrt{\lambda_k}f_k\}$ form an orthonormal basis for $\HHH(F)$.
\end{enumerate}
\end{lem}
\begin{proof}
By definition of $L$,
\[
\la \iota f_k,\iota f_l\ra_\HHH=\la Lf_k,f_l\ra_{F^2}=\lambda_k\delta_{k,l}.
\]
Thus $\{\lambda_k^{-\frac{1}{2}}\iota f_k\}$ form an orthonormal set of $\HHH$. Also, since $h$ is cyclic for $\Tbf$, $\iota$ has dense range in $\HHH$. From this we see that (1) holds.

By definition, $\{f_k\}$ is an orthonormal set in $F^2$. From the cyclicity of $h$,
\[
\III=\{f\in F^2~:~f(\Tbf)=0\}=\{f\in F^2~:~f(\Tbf)h=0\}=\ker L.
\]
Thus $\III^\perp=\overline{\mathrm{ran}L}=\overline{\mathrm{span}\{f_k\}}$. Therefore (2) holds. 

It remains to prove (3). First, by Lemma \ref{lem: f in RKHS criterion} and the expression
\[
F(z,w)=\sum_k\lambda_kf_k(z)\overline{f_k(w)},
\]
we see that each $f_k$ is in $\HHH(F)$, with $\|f_k\|_{\HHH(F)}\leq\lambda_k^{-1/2}$. On the other hand, note that for $c_1,\cdots,c_m\in\CC,~z_1,\cdots,z_m\in\CC^d$,
\[
\sum_{i,j=1}^m c_i\overline{c_j}F(z_j,z_i)=\la Lg,g\ra_{F^2},\quad\text{and}\quad \sum_{i,j=1}^mc_i\overline{c_j}f_k(z_j)\overline{f_k(z_i)}=\left|\la f_k,g\ra_{F^2}\right|^2,\quad\text{where}\quad g=\sum_{i=1}^k c_iK_{z_i}^{F^2}.
\]
For any $c<\lambda_k^{-1/2}$, we may choose $g$ sufficiently close to $f_k$, so that
\[
\sum_{i,j=1}^mc_i\overline{c_j}\left(c^2F(z_j,z_i)-f_k(z_j)\overline{f_k(z_i)}\right)=c^2\la Lg,g\ra_{F^2}-\left|\la f_k,g\ra_{F^2}\right|^2<0.
\]
Consequently, $\|f_k\|_{\HHH(F)}=\lambda_k^{-1/2}, \forall k$. Thus by Lemma \ref{lem: ONB in RKHS}, $\{\sqrt{\lambda_k}f_k\}$ form an orthonormal basis for $\HHH(F)$. This completes the proof.
\end{proof}

\begin{lem}\label{lem: H(F)}
Let $(\Tbf,h)$ be a cyclic commuting $d$-tuple and let $F, \HHH(F)$ be as in Definition \ref{defn: iota L the sequence and F}. Then $\HHH(F)\subset F^2$ as sets of functions. In particular, $\HHH(F)$ consists of analytic functions on $\CC^d$.
\end{lem}
\begin{proof}
By Lemma \ref{lem: ONBs} (3), every function in $\HHH(F)$ has the form 
\[
f=\sum_ka_k\sqrt{\lambda_k}f_k,\quad\text{and}\quad\sum_k|a_k|^2<\infty.
\]
Since $\lambda_k\to\infty$, the series also converges in $F^2$. Since evaluation is continuous in both spaces, the series represents the same function in $\HHH(F)$ and $F^2$. This completes the proof.
\end{proof}

\begin{lem}\label{lem: f_j closed under polynomial}
For any $k$, $f_k\polyzd\subset F^2$.
\end{lem}
\begin{proof}
For any polynomial $p$, any $k, i$,
\[
|\la f_k,\partial_ip\ra_{F^2}|=\lambda_k^{-1/2}|\la L^{1/2}f_k,\partial_ip\ra_{F^2}|=\lambda_k^{-1/2}|\la f_k,L^{1/2}\partial_i p\ra_{F^2}|\leq\lambda_k^{-1/2}\|\partial_ip(\Tbf)h\|\leq C_1\|\partial_ip\|_{r\DD^d}\leq C_2\|p\|_{F^2},
\]
where $C_1, C_2$ are constants independent of $p$. Thus $z_i f_k\in F^2,~\forall k, \forall i$. A similar argument shows that $z^\alpha f_k\in F^2, \forall\alpha\in\NN_0^d$. This completes the proof.
\end{proof}

\begin{proof}[{\bf Proof of Theorem \ref{thm: intro H(F) and T}}]
Statements (1) and (2) are proved in Lemma \ref{lem: H(F)} and Lemma \ref{lem: ONBs} (3), respectively. Next, we prove (3). Note that by Lemma \ref{lem: ONBs} (3) we have $\HHH(F)=L^{1/2}(F^2)=\{\sum_ka_k\sqrt{\lambda_k}f_k~:~\sum_k|a_k|^2<\infty\}$. For any $f, g\in F^2$,
\begin{flalign}\label{eqn: temp 1}
\left|\la L^{1/2}f,z_ig\ra_{F^2}\right|=&\left|\la f, L^{1/2}z_ig\ra_{F^2}\right|\leq\|f\|_{F^2}\|L^{1/2}z_ig\|_{F^2}\nonumber\\
=&\|f\|_{F^2}\|T_ig(\Tbf)h\|_\HHH\leq\|T_i\|\|f\|_{F^2}\|g(\Tbf)h\|_\HHH\nonumber\\
=&\|T_i\|\|f\|_{F^2}\|L^{1/2}g\|_{F^2}.
\end{flalign}
Therefore 
\[
\left|\la L^{1/2}f,z_ig\ra_{F^2}\right|\leq=\|T_i\|\|f\|_{F^2}\|L^{1/2}g\|_{F^2}\leq\|T_i\|\|f\|_{F^2}\|L^{1/2}\|\|g\|_{F^2}.
\]
In other words, the linear functional $\la L^{1/2}f,z_i\cdot\ra_{F^2}$ is bounded in $F^2$. Therefore $\partial_iL^{1/2}f\in F^2$ for any $f\in F^2$.
If $g\in\ker L$, then Equation \eqref{eqn: temp 1} gives $\la\partial_i L^{1/2}f,g\ra_{F^2}=0$. Thus $\partial_i L^{1/2}f\in\overline{\mathrm{ran L}}=\overline{\mathrm{span}}\{f_k\}$. Assume $\partial_i L^{1/2}f=\sum_k b_kf_k$.
Taking $g=\sum_{k=1}^na_kf_k$ in Equation \eqref{eqn: temp 1} gives
\[
|\sum_{k=1}^n\overline{a_k}b_k|=\left|\la \partial_i L^{1/2}f,g\ra_{F^2}\right|\leq C\|L^{1/2}g\|_{F^2}=C\sqrt{\sum_{k=1}^n\lambda_k |a_k|^2},
\]
where $C=\|T_i\|\|f\|_{F^2}$. By duality we have 
\[
\left\|\partial_i L^{1/2}f\right\|_{\HHH(F)}^2=\sum_k|b_k|^2\lambda_k^{-1}\leq C.
\]
Thus we have shown 
\[
\partial_i L^{1/2}f\in \HHH(F),\quad\forall f\in F^2.
\]
In other words, $\partial_i\HHH(F)\subset\HHH(F)$. 
By the closed graph theorem, $\partial_i$ defines a bounded operator on $\HHH(F)$. This proves (3).

For any $k, l$,
\begin{flalign*}
\la T_i^\ast\iota f_k,\iota f_l\ra_\HHH=&\la \iota f_k, T_i\iota f_l\ra_\HHH=\la\iota f_k,\iota z_if_l\ra_\HHH=\la Lf_k,z_if_l\ra_{F^2}\\
=&\la \lambda_kf_k,z_if_l\ra_{F^2}=\la \lambda_k\partial_if_k,f_l\ra_{F^2}=\la\lambda_k\partial_if_k,\lambda_lf_l\ra_{\HHH(F)}\\
=&\la \partial_iU\iota f_k,U\iota f_l\ra_{\HHH(F)}.
\end{flalign*}
This proves (4). The proof is complete.
\end{proof}

\begin{proof}[{\bf Proof of Corollary \ref{cor: intro joint eigenvalue for T*}}]
By Theorem \ref{thm: intro H(F) and T} (4), the joint eigenvalues of $\Tbf^\ast$ equals the joint eigenvalues of $\lambda=(\partial_1,\cdots,\partial_d)$ on $\HHH(F)$. Suppose $f\in\HHH(F)$ is a joint eigenvector for $\mathbf{\partial}$, with eigenvalues $\mathbf{\lambda}=(\lambda_1,\cdots,\lambda_d)$. Then
\[
\partial_if=\lambda_if,\quad\forall i=1,\cdots,d.
\]
Thus $f(z)=ce^{\la z,\bar{\lambda}\ra}$ for some constant $c\in\CC$. In other words, $\lambda\in\sigma_p(\Tbf^\ast)$ if and only if $e^{\la z,\bar{\lambda}\ra}\in\HHH(F)$. Take an orthonormal basis $\{g_l\}$ of $\III$ in $F^2$. Then $\{f_k\}\cup\{g_l\}$ form an orthonormal basis for $F^2$. Therefore
\[
e^{\la z,\bar{\lambda}\ra}=K^{F^2}_{，\bar{\lambda}}(z)=\sum_k\overline{f_k(\bar{\lambda})}f_k(z)+\sum_l\overline{g_l(\bar{\lambda})}g_l(z).
\]
By Lemma \ref{lem: ONBs}, 
\[
e^{\la z,\bar{\lambda}\ra}\in\HHH(F)\quad\Leftrightarrow\quad \sum_k\lambda_k^{-1}|f_k(\bar{\lambda})|^2<\infty,\text{ and }g_l(\bar{\lambda})=0,~\forall l\quad\Leftrightarrow\quad \sum_k\lambda_k^{-1}|f_k(\bar{\lambda})|^2<\infty,\text{ and }\bar{\lambda}\in Z(\III).
\]
This completes the proof.
\end{proof}

\section{A Paley-Weiner-Schwartz Type Theorem}\label{sec: PWS type theorem}
In this section, we give a few examples whose Agler's linear functional are given by distributions. Then we give the proof of Theorem \ref{thm: intro PWS}.
We start with a very simple example.
\begin{exam}[Jordan block]\label{exam: Jordan block}
Let $T$ be the $m\times m$ Jordan block 
\[
T=\begin{bmatrix}    0&0&0&\cdots&0\\1&0&0&\cdots&0\\0&1&0&\cdots&0\\\vdots&\vdots&\vdots&\ddots&\vdots\\0&0&0&\cdots&0
\end{bmatrix}_{m\times m}
\]
viewed as an operator on $\CC^m$, and let $h=e_1$, where $e_{i,j}=\delta_{i,j}$. By direct computation,
\[
\Lambda_{T,h}(z^k\bz^l)=\begin{cases}
\delta_{k,l},&\text{if }k, l\leq m-1\\0,&\text{if }k\text{ or }l\geq m.
\end{cases}
\]
Then it is easy to verify that
\[
\Lambda_{T,h}=\left(I+\partial\bpartial+\frac{1}{4}\partial^2\bpartial^2+\cdots+\frac{1}{\left((m-1)!\right)^2}\partial^{m-1}\bpartial^{m-1}\right)\delta_0=\frac{1}{4^{m-1}\left((m-1)!\right)^2}(1-|z|^2)^{-1}\Delta^{m-1}\delta_0,
\]
where $\delta_0$ denotes the point mass at $0\in\CC$, and the derivation and multiplication are in the sense of distributions. Using properties of Fourier-Laplace transform, we also have
\[
F_{T,h}(z,w)=1+z\bw+\frac{1}{4}z^2\bw^2+\cdots+\frac{1}{\left((m-1)!\right)^2}z^{m-1}\bw^{m-1}.
\]
\end{exam}
Another example comes from the Drury-Arveson space.
\begin{prop}\label{prop: DA space inner product}
There is a unique distribution $u$ supported on the unit sphere $\partial\BB_d$, such that 
\[
\la p, q\ra_{H_d^2}=u(p\bar{q}),\quad\forall p, q\in\polyzd.
\]
Equivalently, let $\Tbf=\Mzbf=(M_{z_1},\cdots,M_{z_d})$ be the tuple of coordinate multiplication operators on $H_d^2$, then $\Lambda_{\Tbf,1}=u$. Indeed,
\[
u=\frac{(-1)^{d-1}}{(d-1)!}\left(R+(d-1)I\right)\cdots(R+I)\sigma,
\]
where $\sigma$ is the normalized surface measure on $\partial\BB_d$, and $R=\sum_i z_i\partial_i$.
\end{prop}

\begin{proof}
Note that $\polytwod$ is dense in $\cinfty(\CC^d)$, the distribution $u$ is unique if it exists. Recall that the Hardy space $H^2(\BB_d)$ is the analytic function space with the same orthogonal basis $\{z^\alpha\}_{\alpha\in\NN_0^d}$, and
\[
\|z^\alpha\|_{H^2(\BB_d)}^2=\frac{\alpha!(d-1)!}{(|\alpha|+d-1)!}.
\]
For polynomials $p, q\in\polyzd$,
\[
\la p, q\ra_{H^2(\BB_d)}=\int_{\partial\BB_d}p\bar{q}~\intd\sigma=\sigma(p\bar{q}),
\]
where in the last expression, $\sigma$ is considered as a bounded linear functional on $C(\partial\BB_d)$. Denote $R=\sum_{i=1}^dz_i\partial_i$, the radial derivative operator. Then $Rz^\alpha=|\alpha|z^\alpha$. 
For any $\alpha, \beta\in\NN_0^d$,
\begin{flalign*}
\la z^\alpha, z^\beta\ra_{H_d^2}=&\delta_{\alpha,\beta}\frac{\alpha!}{|\alpha|!}=\frac{(|\alpha|+d-1)!}{|\alpha|!(d-1)!}\delta_{\alpha,\beta}\frac{\alpha!(d-1)!}{(|\alpha|+d-1)!}=\frac{(|\alpha|+d-1)!}{|\alpha|!(d-1)!}\la z^\alpha, z^\beta\ra_{H^2(\BB_d)}\\
=&\frac{1}{(d-1)!}\la\left(R+(d-1)I\right)\cdots(R+I)z^\alpha, z^\beta\ra_{H^2(\BB_d)}\\
=&\frac{1}{(d-1)!}\sigma\left(\left(R+(d-1)I\right)\cdots(R+I)(z^\alpha\bz^\beta)\right).
\end{flalign*}
It follows that for any $p, q\in\polyzd$, $\la p, q\ra_{H_d^2}=u(p\bar{q})$, where $u$ is the distribution defined by
\[
u(\phi)=\frac{1}{(d-1)!}\sigma\left(\left(R+(d-1)I\right)\cdots(R+I)(\phi)\right),\quad\forall\phi\in\polytwod.
\]
For any distribution $v$ and any $\phi\in\cinfty(\CC^d)$,
\[
v(R\phi)=\sum_{i=1}^dv(z_i\partial_i\phi)=\sum_{i=1}^d(z_iv)(\partial_i\phi)=-\sum_{i=1}^d(\partial_iz_iv)(\phi)=-\sum_{i=1}^d(v+z_i\partial_iv)(\phi)=-((R+dI)v)(\phi).
\]
Thus we may alternatively write
\[
u=\frac{(-1)^{d-1}}{(d-1)!}\left(R+(d-1)I\right)\cdots(R+I)\sigma.
\]
From this we see that $u$ is a distribution supported on $\partial\BB_d$. This completes the proof.
\end{proof}

\begin{rem}\label{rem: subjordan}
In \cite{Helton71}\cite{Helton72}\cite{Helton70}, William Helton initiated research on a class of operators that was later referred to as real Jordan operators. An operator $T$ is real Jordan if it is a sum $T=N+J$, where $N$ is self-adjoint, $J$ is nilpotent, and $NJ=JN$. A real subjordan operator is the restriction of a real Jordan operator on an invariant subspace. In a sequence of works, \cite{Helton70}\cite{Helton71}\cite{Helton72}\cite{BaHe80}\cite{Agler80}, William Helton, Joseph Ball, and Jim Agler gave algebraic characterizations of real Jordan and subjordan operators. In \cite{BaFa90}\cite{BaFa95}, Joseph Ball and Thomas Fanney considered complex Jordan and subjordan operators. Here the condition that $N$ is self-adjoint is replaced by that $N$ is normal. From the function model \cite[Theorem 3.2]{BaFa95}, it is easy to see that for any complex Jordan operator $T$ on $\HHH$ and any $h\in\HHH$, $\Lambda_{T,h}$ is a distribution of certain form. 
Later, in \cite{Agler90}\cite{AgSt95I}\cite{AgSt95II}\cite{AgSt96III}, Jim Agler and Mark Stankus initiated research on a parallel theory called $m$-isometries. This has become an active area of research even today. There are a lot of references on this topic. Here we list some of them as examples, \cite{BeMaNe10}\cite{Richter91}\cite{GlRi06}\cite{GuSt15}\cite{BaSu19}\cite{Le20}.
In particular, in \cite{GlRi06}, Jim Gleason and Stefan Richter developed the theory of multi-variable $m$-isometries and showed that $n$-shift is an $n$-isometry.
\end{rem}

In the rest of this section, we give the proof of Theorem \ref{thm: intro PWS}.
First, let us recall some basic facts in distribution theory. Our main source of references is \cite{HormanderPDEI2003}. Recall the Fr\'{e}chet space
\[
\cinfty(\RR^d):=\{f: \RR^d\to\CC\text{ smooth}\}
\]
is equipped with the semi-norms
\[
\|f\|_{N,K}:=\max\{|D^\alpha f(x)|~:~|\alpha|\leq N, x\in K\},\quad\text{where }N\in\NN_0^d, K\subset\RR^d\text{ is compact}.
\]
The set of all compactly supported distributions on $\RR^d$ coincides with the space of continuous linear functionals on $\cinfty(\RR^d)$. 
Suppose $u$ is a compactly supported distribution on $\RR^d$. Recall that its Fourier-Laplace transform is the entire analytic function on $\CC^d$, defined by
\[
\hat{u}(\zeta)=u_x\left(e^{-i\la\zeta,x\ra}\right),\quad\forall\zeta\in\CC^d.
\]
\begin{defn}[\cite{HormanderPDEI2003} Definition 4.3.1]\label{defn: supporting function real}
Suppose $E$ is a compact set in $\RR^d$. Define
\[
H_E(\xi)=\sup_{x\in E}\la x,\xi\ra,\quad\xi\in\RR^d.
\]
The function $H_E$ is called the \emph{supporting function} of $E$.
\end{defn}
The Paley-Wiener-Schwartz Theorem characterizes compactly supported distributions through their Fourier-Laplace transforms. 
\begin{thm}[Paley-Wiener-Schwartz \cite{HormanderPDEI2003} Theorem 7.3.1]\label{thm: PWS thm}
Let $E$ be a convex compact subset of $\RR^d$ with supporting function $H$. If $u$ is a distribution of order $N$ with support contained in $E$, then
\begin{equation}\label{eqn: PWS real}
|\hat{u}(\zeta)|\leq C(1+|\zeta|)^Ne^{H(\mathrm{Im}\zeta)},\quad\zeta\in\CC^d.
\end{equation}
Conversely, every entire analytic function in $\CC^d$ satisfying an estimate of the form \eqref{eqn: PWS real} is the Fourier-Laplace transform of a distribution with support contained in $E$.
\end{thm}

Suppose $u$ is a compactly supported distribution on $\CC^d=\RR^{2d}$. Then $\hat{u}$ is an entire function on $\CC^{2d}$. The definition of $\hat{u}$ does not reflect the complex structure of the underlying space $\CC^d$. Thus it makes sense to define
\begin{equation}\label{eqn: Fu definition}
F_u(z,w)=u_{\lambda}\left(e^{\la\lambda,w\ra+\la z,\lambda\ra}\right),\quad\forall z, w\in\CC^d.
\end{equation}
It is straightforward to verify the following identities
\begin{equation}\label{eqn: identity btw uhat and Fu}
\hat{u}(\zeta_1,\zeta_2)=F_u(\frac{\zeta_2-i\zeta_1}{2},\frac{i\bar{\zeta_1}-\bar{\zeta_2}}{2}),\quad\text{and}\quad F_u(z,w)=\hat{u}\left((i(\bw+z),-\bw+z)\right),
\end{equation}
for any $\zeta_1,\zeta_2,z,w\in\CC^d$.

\begin{defn}\label{defn: complex supporting function}
For a compact convex set $K$ in $\CC^d$, we define its \emph{(complex) supporting function} $H_K$ to be
\[
H_K(z)=\sup_{\lambda\in K}\mathrm{Re}\la\lambda,z\ra,\quad z\in\CC^d.
\]
It is the complex analogue of the real supporting function in Definition \ref{defn: supporting function real}.
\end{defn}

\begin{proof}[{\bf Proof of Theorem \ref{thm: intro PWS}}]
Suppose $\Lambda_{\Tbf,h}=u$ is a distribution of order $N$, supported in $K$. Identify $\CC^d$ with $\RR^{2d}$, and let 
\[
E=\{(x,y)~:~x, y\in\RR^d, x+iy\in K\}.
\] 
Denote $H_E$ the supporting function of $E$ defined as in Definition \ref{defn: supporting function real}. By Theorem \ref{thm: PWS thm},
\[
\left|\hat{u}(\zeta)\right|\leq C(1+|\zeta|)^N e^{H_E(Im\zeta)},\quad\forall\zeta\in\CC^{2n}.
\]
Thus for some $C_1>0$,
\begin{flalign*}
|F_{\Tbf,h}(z,w)|=&|F_u(z,w)|=|\hat{u}\left((i(\bw+z),-\bw+z)\right)|\\
\leq&C_1(1+|z|+|w|)^Ne^{H_E((\mathrm{Re}(\bar{w}+z),\mathrm{Im}(-\bar{w}+z)))}\\
=&C_1(1+|z|+|w|)^Ne^{H_E((\mathrm{Re}(w+z),\mathrm{Im}(w+z)))}.
\end{flalign*}
Also,
\begin{flalign*}
H_E((\mathrm{Re}(w+z),\mathrm{Im}(w+z)))
=&\sup_{x\in E}\left(\la x_1,\mathrm{Re}(w+z)\ra+\la x_2, \mathrm{Im}(w+z)\ra\right)\\
=&\sup_{\lambda\in K}\mathrm{Re}\la\lambda, w+z\ra\\
=&H_K(z+w).
\end{flalign*}
Combining the two inequlities above gives Inequality \eqref{eqn: intro F PWS}.

Now suppose $F_{\Tbf,h}$ satisfies Inequality \eqref{eqn: intro F PWS}. Define $g: \CC^{2d}\to\CC$,
\[
g(\zeta_1,\zeta_2)=F_{\Tbf,h}\left(\frac{\zeta_2-i\zeta_1}{2},\frac{i\bar{\zeta_1}-\bar{\zeta_2}}{2}\right).
\]
Then, similar as above, $g$ is an entire function satisfying the estimate
\[
|g(\zeta)|\leq C_2(1+|\zeta|)^N e^{H_E(Im\zeta)},\quad\forall\zeta\in\CC^{2d}
\]
for some $C_2>0$. By Theorem \ref{thm: PWS thm}, there is some distribution $u$ of order $N$, supported in $E$ such that $g=\hat{u}$. By \eqref{eqn: identity btw uhat and Fu}, $F_{\Tbf,h}=F_u$. From this it is easy to see that $\Lambda_{\Tbf,h}=u$. This completes the proof.
\end{proof}

\section{The Matrix Case}\label{sec: the matrix case}
In the case when $T$ consists of matrices, the characterization in Theorem \ref{thm: intro PWS} can be more explicit. 
We begin by fixing a few terminologies. 
\begin{defn}\label{defn: Jordan}
\begin{enumerate}
    \item We say that a commuting tuple $\Tbf=(T_1,\cdots,T_d)$ on $\HHH$ is \emph{nilpotent} if for some $\mbf\in\NN_0^d$, $T_i^{m_i}=0, i=1,\cdots,d$. It is easy to see that a cyclic nilpotent tuple must consist of matrices.
    \item We say that a commuting tuple of matrices $\Jbf=(J_1,\cdots,J_d)$ is a \emph{tuple of (generalized) Jordan blocks} at $\lambda\in\CC^d$ if for some commuting tuple $\Nbf$ of nilpotent matrices, $\Jbf=\lambda\Ibf+\Nbf$, i.e., $J_i=\lambda_i I+N_i, i=1,\cdots,d$. 
    \item We say that a commuting tuple $T=(T_1,\cdots,T_d)$ of matrices is \emph{a (generalized) Jordan tuple} if there is a unitary operator $U$ such that the tuple $U^\ast TU:=(U^\ast T_1 U,\cdots, U^\ast T_d U)$ is a direct sum of tuples of (generalized) Jordan blocks.
\end{enumerate}
\end{defn}

In this section, we prove the following theorem.
\begin{thm}\label{thm: matrix case}
Suppose $(\Tbf,h)$ is a cyclic commuting tuple on $\CC^m$. Then the following are equivalent.
\begin{enumerate}
    \item $\Lambda_{\Tbf,h}$ is a distribution.
    \item $\Tbf$ is a Jordan tuple.
\end{enumerate}
In this case, $\Lambda_{\Tbf,h}$ is supported at $\sigma(\Tbf)$, the joint Taylor spectrum of $\Tbf$, which consists of finitely many points.
\end{thm}

To present the proof, we find it convenient to treat a linear functional on $\polytwod$ as a generalized class of distribution, and extend some operations on distributions to this generalized class.

\begin{defn}\label{defn: linear functionals as distributions}
Suppose $\Lambda$ is a linear functional on $\polytwod$, $\phi\in\polytwod$, $\lambda\in\CC^d$, and $\alpha, \beta\in\NN_0^d$.
\begin{enumerate}
    \item Define the linear functional $\phi\Lambda$ by
    \[
    \phi\Lambda: \polytwod\to\CC,\quad\phi\Lambda(\psi)=\Lambda(\phi\psi).
    \]
    \item Define the translation $\tau_\lambda\Lambda$ to be the linear functional
    \[
    \tau_\lambda\Lambda: \polytwod\to\CC,\quad\tau_\lambda\Lambda(\psi)=\Lambda(\tau_{-\lambda}\psi),
    \]
    where $\tau_\lambda\psi(z):=\psi(z-\lambda)$.
    \item Define the differentiation $\partial^\alpha\bpartial^\beta\Lambda$ to be the linear functional
    \[
    \partial^\alpha\bpartial^\beta\Lambda: \polytwod\to\CC,\quad\partial^\alpha\bpartial^\beta\Lambda(\psi)=(-1)^{|\alpha|+|\beta|}\Lambda(\partial^\alpha\bpartial^\beta\psi).
    \]
\end{enumerate}
\end{defn}

The following lemmas are straightforward to check.

\begin{lem}\label{lem: LambdaTh for direct sums}
Suppose $\Tbf$ and $\Sbf$ are commuting $d$-tuples on $\HHH$, $\EEE$, respectively, and $h\in\HHH, e\in\EEE$. Then
\[
\Lambda_{\Tbf\oplus\Sbf, h\oplus e}=\Lambda_{\Tbf,h}+\Lambda_{\Sbf,e}.
\]
\end{lem}
We remark here that $h\oplus e$ may no longer be cyclic for $\Tbf\oplus\Sbf$. However, the definition of $\Lambda_{\Tbf,h}$ carries over to the non-cyclic case with no difference.
\begin{lem}\label{lem: LambdaTh for translation}
Suppose $(\Tbf, h)$ is a commuting $d$-tuple on $\HHH$ and $h\in\HHH$. For $\lambda\in\CC^d$, denote
\[
\Tbf+\lambda\Ibf:= (T_1+\lambda_1 I,\cdots, T_d+\lambda_dI).
\]
Then $\Lambda_{\Tbf+\lambda\Ibf, h}=\tau_\lambda\Lambda_{\Tbf,h}$.
\end{lem}

\begin{lem}\label{lem: LambdaTh for p(T)h}
Suppose $\Tbf$ is a commuting $d$-tuple on $\HHH$, $h\in\HHH$, and $p\in\polyzd$. Then
\[
\Lambda_{\Tbf,p(\Tbf)h}=|p|^2\Lambda_{\Tbf,h}.
\]
\end{lem}

\begin{lem}\label{lem: p(T)=0 equivalent}
Suppose $(\Tbf,h)$ is a cyclic commuting tuple. Then for any $p\in\polyzd$,
\[
p(\Tbf)=0\quad\Leftrightarrow\quad p\Lambda_{\Tbf,h}=0,
\]
and
\[
p(\Tbf)\text{ is self adjoint}\quad\Leftrightarrow\quad(p-\bar{p})\Lambda_{\Tbf,h}=0.
\]
\end{lem}

We start by proving a special case of Theorem \ref{thm: matrix case}.

\begin{lem}\label{lem: LambdaTh for nilpotent}
Suppose $(\Tbf,h)$ is a cyclic commuting $d$-tuple on $\HHH$. Then the following are equivalent.
\begin{enumerate}
    \item $\Tbf$ is nilpotent;
    \item $\Lambda_{\Tbf,h}$ is a distribution supported at $0$;
    \item there exists a finite collection $\{q_k\}\subset\polyzd$ such that
\[
\Lambda_{\Tbf,h}=\sum_k q_k(\partial)\overline{q_k(\partial)}\delta_0;
\]
\item there exists a finite collection $\{p_k\}\subset\polyzd$ such that 
\[
L_{\Tbf,h}=\sum_kp_k\otimes p_k.
\]
\end{enumerate}
\end{lem}
Here for $f, g\in F^2$, $f\otimes g$ denotes the operator
\[
f\otimes g(r)=\la r, g\ra_{F^2}f,\quad\forall r\in F^2.
\]
\begin{proof}
The fact that (3)$\Leftrightarrow$(4) is obvious once one notices
\[
\la f, p\ra_{F^2}=\left(p^\ast(\partial)f\right)(0),\quad\forall p\in\polyzd, \forall f\in F^2,\quad\text{where }p^\ast(z)=\sum_\alpha\overline{c_\alpha}z^\alpha\text{ for }p(z)=\sum_\alpha c_\alpha z^\alpha.
\]
(1)$\Rightarrow$(2): Assume that for $\mbf\in\NN_0^d, T_i^{m_i}=0, i=1,\cdots, d$. Recall in Definition \ref{defn: iota L the sequence and F} and the GNS-type construction below, we defined the semi-inner product
\[
\la p, q\ra_{\Tbf,h}=\la p(\Tbf)h,q(\Tbf)h\ra_\HHH,\quad\forall p, q\in\polyzd,
\]
and the ideal
\[
I:=\{p\in\polyzd~:~\|p\|_{\Tbf,h}=0\}=\{p\in\polyzd~:~p(\Tbf)=0\}.
\]
Still denote by $\la\cdot,\cdot\ra_{\Tbf,h}$ the induced inner product on
induces an inner product on $\polyzd/I$. By assumption, $z_i^{m_i}\in I$. So $I$ is a primary ideal with a one point zero locus $\{0\}$, and $\polyzd/I$ is finite-dimensional. Since every positive definite matrix has a decomposition into rank one positive definite matrices, there exist a finite collection of linear functionals $\{\ell_k\}$ on $\polyzd$, each vanishing on $I$, such that
\[
\la p, q\ra_{\Tbf,h}=\sum_k\ell_k(p)\overline{\ell_k(q)},\quad\forall p, q\in\polyzd.
\]
It suffices to characterize all such linear functionals. Here we use the theory of characteristic spaces defined by Kunyu Guo (cf. \cite{ChGuAnalyticHilbertModules}\cite{Guo99}\cite{Guo2000}). As in \cite[Section 2.1]{ChGuAnalyticHilbertModules}, define the characteristic space of $I$ at $0$ to be
\[
I_0=\{q\in\polyzd~|~\left(q(\partial)p\right)(0)=0, \forall p\in I\}.
\]
By \cite[Corollary 2.1.2]{ChGuAnalyticHilbertModules},
\[
I=\{p\in\polyzd~:~\left(q(\partial)p\right)(0)=0, \forall q\in I_0\}.
\]
Thus for any linear functional $\ell$ on $\polyzd$ that vanishes on $I$, there is $q\in I_0$ such that
\[
\ell(p)=\left(q(\partial)p\right)(0),\quad\forall p\in\polyzd.
\]
Therefore for some finite subset $\{q_k\}\subset I_0$, for all $p, q\in\polyzd$,
\[
\Lambda_{\Tbf,h}(p\bar{q})=\la p, q\ra_{\Tbf,h}=\sum_k\left(q_k(\partial)p\right)(0)\overline{\left(q_k(\partial )q\right)(0)}=\sum_{k}\left(q_k(\partial)q_k^\ast(\bpartial)(p\bar{q})\right)(0)=\left(\sum_k\tilde{q_k}(\partial)\tilde{q_k}^\ast(\bpartial)\delta_0\right)(p\bar{q}).
\]
Here $\tilde{q}(z)=\sum_\alpha(-1)^{|\alpha|}c_\alpha z^\alpha$ for $q(z)=\sum_\alpha c_\alpha z^\alpha$. This proves (2).

(2)$\Rightarrow$(3): If $\Lambda_{\Tbf,h}$ is a distribution supported at $0$. By \cite[Theorem 2.3.4]{HormanderPDEI2003}, $\la\cdot,\cdot\ra_{\Tbf,h}$ vanishes in a primary ideal with single point zero locus $\{0\}$. Repeating the above argument gives (3).

(4)$\Rightarrow$(1): Let $m$ be an integer that is larger than the degree of all $p_k$. Then it is straightforward to check that $T_i^m=0$ for all $i=1,\cdots,d$. This completes the proof.
\end{proof}

The following lemma follows directly from Lemma \ref{lem: LambdaTh for translation} and Lemma \ref{lem: LambdaTh for nilpotent}.
\begin{lem}\label{lem: LambdaTh for Jordan block}
Suppose $(\Tbf,h)$ is a cyclic commuting $d$-tuple on $\HHH$ and $\lambda\in\CC^n$. Then the following are equivalent.
\begin{enumerate}
    \item $\Tbf$ is a Jordan block at $\lambda$;
    \item $\Lambda_{\Tbf,h}$ is a distribution supported at $\lambda$;
    \item there exists a finite collection $\{q_k\}\subset\polyzd$ such that
\[
\Lambda_{\Tbf,h}=\sum_k q_k(\partial)\overline{q_k(\partial)}\delta_\lambda;
\]
\item there exists a finite collection $\{p_k\}\subset\polyzd$ such that 
\[
L_{\Tbf,h}=\sum_kp_kK_\lambda\otimes p_kK_\lambda.
\]
\end{enumerate}
\end{lem}

With the preparations above, we are ready to prove Theorem \ref{thm: matrix case}.

\begin{proof}[{\bf Proof of Theorem \ref{thm: matrix case}}]
The fact that (2) implies (1) follows from Lemma \ref{lem: LambdaTh for Jordan block} and Lemma \ref{lem: LambdaTh for direct sums}.
We show that (1) implies (2). Suppose $\Lambda_{\Tbf,h}$ is a distribution. As in the proof of Lemma \ref{lem: LambdaTh for nilpotent}, let 
\[
I=\{p\in\polyzd~:~p(\Tbf)=0\}=\{p\in\polyzd~:~\|p\|_{\Tbf,h}=0\}.
\]
By Lemma \ref{lem: p(T)=0 equivalent},
\[
p\Lambda_{\Tbf,h}=0,\quad\forall p\in I.
\]
From this we see that $\Lambda_{\Tbf,h}$ is supported in the zero locus of $I$, which equals $\sigma(\Tbf)$.

Let $N$ be the order of the distribution $\Lambda_{\Tbf,h}$. Write $\sigma(\Tbf)=\{\lambda_k\}_{k=1}^m$. For each $k$, choose a polynomial $p_k\in\polyzd$ such that
\[
p_k(\lambda_l)=\delta_{k,l},\quad\text{and}\quad\partial^\alpha p_k(\lambda_l)=0, \forall k, l, 0<|\alpha|\leq N.
\]
Then it is easy to verify that
\[
\left(p_k^2-p_k\right)\Lambda_{\Tbf,h}=0,\quad\left(p_k-\overline{p_k}\right)\Lambda_{\Tbf,h}=0,\quad p_kp_l\Lambda_{\Tbf,h}=0, \forall k\neq l,\quad \left(\sum_k p_k\right)\Lambda_{\Tbf,h}=\Lambda_{\Tbf,h}.
\]
Again, by Lemma \ref{lem: p(T)=0 equivalent},
\[
p_k(\Tbf)^2=p_k(\Tbf),\quad p_k(\Tbf)=\left(p_k(\Tbf)\right)^\ast,\quad p_k(\Tbf)p_l(\Tbf)=0, \forall k\neq l,\quad\sum_kp_k(\Tbf)=I.
\]
In other words, the operators $p_k(\Tbf)$ are pairwise orthogonal projections, and $\sum_k p_k(\Tbf)=I$. Since each $p_k(\Tbf)$ commutes with each $T_i$, its range space is a joint reducing subspace of $\Tbf$. Denote $\HHH_k$ the range of $p_k(\Tbf)$ and $T_{i,k}$ the compression of $T_i$ to $\HHH_k$, $\Tbf_k=(T_{1,k},\cdots,T_{d,k})$. Then $\Tbf=\oplus_k \Tbf_k$. 

We show that each $\Tbf_k$ is a tuple of Jordan blocks. Let $h_k=p_k(\Tbf)h\in\HHH_k$. Clearly $h_k$ is cyclic for $\Tbf_k$. By Lemma \ref{lem: LambdaTh for p(T)h},
\[
\Lambda_{\Tbf_k,h_k}=\Lambda_{\Tbf,h_k}=\Lambda_{\Tbf,p_k(\Tbf)h}=|p_k|^2\Lambda_{\Tbf,h}.
\]
Thus $\Lambda_{\Tbf_k,h_k}=|p_k|^2\Lambda_{\Tbf,h}$ is a distribution supported at a single point $\{\lambda_k\}$. By Lemma \ref{lem: LambdaTh for Jordan block}, $\Tbf_k$ is a tuple of Jordan blocks. Therefore $\Tbf=\oplus_k\Tbf_k$ is Jordan. This completes the proof.
\end{proof}

\begin{rem}
The question of whether $\Lambda_{\Tbf,h}$ is a distribution may be related to the theory of generalized scalar operators, developed by Foias (cf. \cite{CoFo68SpectralOp}\cite{EsPu96Spectral}\cite{LaNe2000LST}). Recall that a generalized scalar operator (cf. \cite[Definition 1.4.9]{LaNe2000LST}) is a bounded linear operator $T$ on a Banach space $X$ for which there exists an algebra homomorphism $\Phi: \cinfty(\CC)\to \BBB(X)$ such that $\Phi(1)=I$ and $\Phi(z)=T$. However, this homomorphism does not necessarily preserve adjoints. From the definition, it is also immediate that this property is stable under similarity equivalence, and in particular, any finite matrix is generalized scalar. Meanwhile, as we have seen in Theorem \ref{thm: matrix case}, the linear functional $\Lambda_{\Tbf,h}$ is only stable under unitary equivalence and may change dramatically under similarity equivalence. To the authors best knowledge, the results in this paper are not covered by the local spectral theory.
\end{rem}

\bibliographystyle{plain}
\bibliography{reference}
	
\end{document}